\documentclass[12pt,amssymb,amscd]{amsart}
\usepackage{color}
\usepackage{amssymb}
\usepackage{amsmath}
\usepackage{amsthm} 
\usepackage{mathrsfs}                                                           
\usepackage{slashbox}  
\usepackage[all]{xy}
\usepackage{graphicx}
\usepackage{wrapfig} 

\def\si{{\Sigma}}

\def\2{{1\over 2}}

\newcommand{\rf}[1]{(\ref{#1})}
\def\b{\bar}

\newcommand{\sw}{\mathcal{W}}
\newcommand{\sz}{\mathcal{Z}}
\newcommand{\cla}{\check{\lambda}}
\newcommand{\crho}{\check{\rho}}
\newcommand{\mf}{\mathcal{F}}
\newcommand{\sop}{s{\rm Op}_{G}}
\newcommand{\smop}{s{\rm MOp}_G}

\renewcommand{\t}{\tilde}
\newcommand{\p}{\partial}

\def\p{\partial}

\def\b{\bar}

\def\<{\langle}
\def\>{\rangle}
\def\+{\dagger}

\begin{document}
\title{Superopers on supercurves}
\author{Anton M. Zeitlin}
\address{ 
\newline Department of Mathematics,\newline
Columbia University,
\newline 2990 Broadway, New York,\newline
NY 10027, USA.
\newline
Max Planck Institute for Mathematics,\newline
Vivatsgasse 7, Bonn, 53111, Germany\newline
IPME RAS, V.O. Bolshoj pr., 61, 199178, 
St. Petersburg \newline
zeitlin@math.columbia.edu\newline
http://math.columbia.edu/$\sim$zeitlin \newline
http://www.ipme.ru/zam.html  }

\begin{abstract}
In this note, we 
introduce the generalization of opers (superopers) for a certain class of superalgebras, 
which have pure odd simple root system. 
We study in detail $SPL_2$-superopers and in particular derive  
the corresponding Bethe ansatz equations, 
which describe the spectrum of $osp(2|1)$ Gaudin model.    
\end{abstract}
\maketitle
\section{Introduction}
Opers are necessary ingredients in the study of the geometric Langlands correspondence 
(see e.g. \cite{FL}). They also play important role in many aspects of mathematical 
physics. For example, opers are very important in the theory of integrable systems, and recently they became a necessary component even in the modern 
Quantum Field Theory approaches to the knot theory (see e.g. \cite{W}).       

Originally, opers were studied locally in the seminal paper of Drinfeld and Sokolov \cite{DS}
as gauge equivalence classes of certain  
differential operators with values in some simple Lie algebra, which are the L-operators of the generalized 
Korteweg-de Vries (KdV) integrable models.  Later, Beilinson and 
Drinfeld generalized  this local object making it coordinate independent \cite{BD}.  
Namely, a $G$-oper on a smooth curve $\Sigma$, where $G$ is a simple algebraic group of the adjoint type with the Lie algebra $\mathfrak{g}$, 
is a triple $(\mathcal{F},\mathcal{F}_B,\nabla$), where 
$\mathcal{F}$ is a $G$-bundle over $\Sigma$, $\mathcal{F}_B$ is its $B$-reduction with respect to 
Borel subgroup $B$, and $\nabla$ is a flat connection, which behaves in 
a certain way with respect to $\mathcal{F}_B$. For example, in the case of $PGL_2$-oper, this 
condition just means that the reduction $\mathcal{F}_B$ is nowhere preserved by this connection.
Moreover, it appears, following the results of Drinfeld and Sokolov, that the space of $G$-opers is 
equivalent to a certain space of scalar pseudodifferential operators. In the $PGL_2$ case, the resulting 
space of scalar operators is just a family of Sturm-Liouville operators and the connection transformation 
properties allows to consider them on all $\Sigma$ as projective connections.

A really interesting story starts when we allow opers to have regular singularities. It turns out that the opers on the projective line can be described via the Bethe ansatz equations for the Gaudin model corresponding to the Langlands dual Lie algebra \cite{Fe}, \cite{Fb}. 
An important object on the way to understand this relation is the so-called Miura oper, which was introduced by E. Frenkel \cite{Fb}. A $Miura$ $oper$ is an oper with one extra constraint: the connection preserves another 
$B$-reduction of $\mathcal{F}$, which we call $\mathcal{F}'_B$.   
The space of the  
Miura opers, associated to a given $G$-oper with trivial monodromy, is isomorphic to the flag manifold $G/B$. 
If the reduction ${\mathcal{F}'_B}$ corresponds to the point in a big cell of  $G/B$, then such a Miura oper is called $generic$. It was shown by E. Frenkel that any Miura oper is generic on the punctured disc and that there is an isomorphism between the space of generic Miura opers on the open neighborhood with certain $H$-bundle connections ($H=B/[B,B]$) \cite{Fb}. 
The map from $H$-connections to $G$-opers is just a generalization of the standard Miura transformation in the theory of KdV integrable models.

By means of the above relation with the $H$-connections, 
it was proved for $PGL_2$-oper in \cite{Fe} and then 
generalized to the higher rank in \cite{FFR}, \cite{Fb} that the eigenvalues of the Gaudin model for a Langlands dual Lie algebra $\mathfrak{g}^L$ can be described by 
the $G$-opers on $\mathbb{C}P^1$ with given regular singularities and trivial monodromy. 
Namely, the consistency conditions for the $H$-connections underlying such opers coincide with the Bethe ansatz equations for the Gaudin Model.

In this article, we are trying to generalize some of the above notions and results on the level of superalgebras. We define an analogue of the oper in the case of 
supergroups which allow the pure fermionic family of simple roots on a super Riemann surface, following some local considerations of \cite{inami}, \cite{DG}, \cite{kz}. We call such objects $superopers$, and in some sense 
they turn out to be ``square roots'' of standard opers. Unfortunately for all other superalgebras, the resulting formalism allows only locally defined objects (on a formal superdisc). 
We study in detail the simplest nontrivial case of superoper, related to the group $SPL_2$ (see e.g. \cite{crane}), related to superprojective transformations, and explicitly establish the relation between $osp(2|1)$ Gaudin model studied in \cite{kulish} and the 
$SPL_2$-oper on super Riemann sphere with given regular singularities. 

In section 2 we explain the relation between super projective structures on super Riemann surface and the supersymmetric version of the Sturm-Liouville operator. Then we relate it to the flat connection 
on $SPL_2$-bundle which will give us the first example of superoper. 

In Section 3 we use this experience to generalize the notion of superoper to the case of higher rank simple supergroups. However, only the supergroups which permit a pure fermionic system of simple roots allow us to construct a globally defined object on a super Riemann surface.
 We define Miura superopers and superopers with regular singularities   
in section 4. There we study the consistency conditions for the superopers on the superconformal sphere and derive the corresponding Bethe equations. 
We compare the results with the $osp(2|1)$ Gaudin model 
and find that the Bethe ansatz equations coincide with the ``body'' part of the consistency condition for corresponing $SPL_2$ Miura superopers. 

Some remarks and open questions are given in section 5.\\

\noindent{\bf Acknowledgments.} I am very grateful to I. Penkov for useful discussions and to D. Leites for pointing out important references. I am indebted to E. Frenkel and E. Vishnyakova for comments on the manuscript.

\section{Superprojective structures, super Sturm-Liouville operator and\\ $osp(2|1)$ superoper}

\noindent {\bf 2.1. Super Riemann surfaces and superconformal transformations.}
We remind that a $supercurve$ of dimension $(1|1)$  (see e.g. \cite{br}) over some Grassman algebra $S$ is a pair $(X,\mathcal{O}_{X})$, where $X$ is a topological space and $\mathcal{O}_X$ is a sheaf of supercommutative $S$-algebras over $X$ such that $(X,\mathcal{O}^{\rm{red}}_{X})$ is an algebraic curve (where $\mathcal{O}^{\rm{red}}_{X}$ is obtained from $\mathcal{O}_X$ by quoting out nilpotents) and for some open sets $U_{\alpha}\subset X$ and some linearly independent elements $\{\theta_{\alpha}\}$ we have $\mathcal{O}_{U_{\alpha}}=\mathcal{O}^{\rm{red}}_{U_{\alpha}}\otimes S[\theta_{\alpha}]$. These open sets $U_{\alpha}$ serve as coordinate neighborhoods for supercurves with coordinates $(z_{\alpha}, \theta_{\alpha})$.  The coordinate transformations on the overlaps $U_{\alpha}\cup U_{\beta}$ 
are given by the following formulas $z_{\alpha}=F_{\alpha\beta}(z_{\beta}, \theta_{\beta})$, $\theta_{\alpha}=\Psi_{\alpha\beta}(z_{\beta}, \theta_{\beta})$, where $F_{\alpha\beta}$, $\Psi_{\alpha\beta}$ are even and odd functions correspondingly. 
A super Riemann surface $\si$ over some Grassmann algebra $S$ 
(for more details see e.g. \cite{Wl}) is a supercurve of dimension $1|1$ over $S$, with one more extra structure: there is 
a subbundle $\mathcal{D}$ of $T\Sigma$ of dimension $0|1$, such that for any nonzero section $D$ of $\mathcal{D}$ on an 
open subset $U$ of $\si$, $D^2$ is nowhere proportional to $\mathcal{D}$, i.e. we have the exact sequence:
\begin{eqnarray}\label{exact}
0\to \mathcal{D}\to T\si\to \mathcal{D}^2\to 0.
\end{eqnarray}
One can pick the holomorphic local coordinates in such a way that this odd vector field 
will have the form $f(z,\theta)D_{\theta}$, where $f(z,\theta)$ is a non vanishing function and:
\begin{eqnarray}
D_{\theta}=\partial_{\theta}+\theta\partial_z, \quad D_{\theta}^2=\partial_z.
\end{eqnarray}
Such coordinates are called $superconformal$. The transformation between two superconformal coordinate systems 
$(z, \theta)$, $(z', \theta')$ is determined by the condition that $\mathcal{D}$ should be preserved,
i.e.:  
\begin{eqnarray}
D_{\theta}=(D_{\theta}\theta') D_{\theta'},
\end{eqnarray}
so that the constraint on the transformation coming from the local change of coordinates is  
$D_{\theta} z'-\theta'D_\theta \theta'=0$.
An important nontrivial example of a super Riemann surface is the Riemann super sphere $SC^*$: there are two  
charts $(z, \theta)$, $(z, \theta')$ so that
\begin{eqnarray}
z'=-\frac{1}{z},\quad \theta'=\frac{\theta}{z}.
\end{eqnarray}
There is a group of superconformal transformations, usually denoted as $SPL_2$ 
which acts transitively on $SC^*$ as follows: 
\begin{eqnarray}\label{transff}
&&z\to \frac{az+b}{cz+d}+\theta\frac{\gamma z+\delta}{(cz+d)^2}, \nonumber\\
&&\theta\to \frac{\gamma z+\delta}{cz+d}+\theta\frac{1+\2 \delta\gamma}{cz+d}, 
\end{eqnarray}
where $a, b, c,d$ are even, $ad-bc=1$, and $\gamma, \delta$ are odd. The Lie algebra of this group is isomorphic to $osp(2|1)$. 

Let us introduce two more notions which we will use in the following. From now on let us call the sections of 
$\mathcal{D}^n$  the $superconformal$ $fields$ of dimension $-n/2$. 
In particular, taking the dual of the exact sequence \ref{exact}, 
we find that a bundle of superconformal fields of dimension 1 (i.e. $\mathcal{D}^{-2}$) 
is a subbundle in $T^*\si$. Considering the superconformal coordinate system, a nonzero section of 
this bundle is generated by $\eta=dz-\theta d\theta$, which is orthogonal to $D_{\theta}$ under standard pairing.

At last, we introduce one more notation. For any 
element $A$ which belongs to some free module over $S[\theta]$, where $\theta$ is a local odd coordinate, we denote the body of this element (i.e. $A$ is stripped of the dependence on the odd variables) as $\bar{A}$.\\

\noindent {\bf 2.2. Superprojective structures and superprojective connections.} 
Let us at first define what a superprojective connection is. We consider the following differential operator, defined locally with coordinates $(z,\theta)$:
\begin{eqnarray}\label{loc}
D^3_{\theta}-\omega(z,\theta).
\end{eqnarray}
The following proposition holds.\\

\noindent {\bf Proposition 2.1.} \cite{mathieu} {\it Formula \rf{loc} defines the operator $L$, such that
\begin{eqnarray}
L: \mathcal{D}^{-1}\to \mathcal{D}^{2}
\end{eqnarray}
iff the transformation of $\omega$ on the overlap of two coordinate charts $(z,\theta)$, $(z', \theta')$
is given by the following expression:
\begin{eqnarray}
\omega(z, \theta)=\omega(z',\theta')(D_{\theta}\theta')^3+\{\theta';z,\theta\} 
\end{eqnarray}
where
\begin{eqnarray}
\{\theta';z,\theta\}=\frac{\partial^2_z\theta'}{D_{\theta}\theta}-
2\frac{\partial_z\theta D^3_{\theta}\theta'}{(D_{\theta}\theta')^2}
\end{eqnarray}
is a supersymmetric generalization of Schwarzian derivative.}
\\

 One can show that the only 
coordinate transformations for which the super Schwarzian derivative vanishes, are the fractional linear transformations 
\rf{transff}.

Let us consider the covering of $\si$ by open subsets, so that the transition functions are given by \rf{transff}. Two such coverings are considered equivalent if their union has the same property of transition functions. The corresponding equivalence classes are called {\it superprojective structures}. 

It appears that like in the pure even case, there is a bijection between super projective connections and super projective 
structures. For a given super projective structure one can define a superprojective connection by assigning operator $D_{\theta}^3$ in every coordinate chart. From Proposition 2.1 we find that the resulting object is defined globally on $\si$. 
On the other hand, given a super projective connection on $\si$, one can consider the following linear problem:
\begin{eqnarray}
(D^3_{\theta}-\omega(z,\theta))\psi(z, \theta)=0.
\end{eqnarray}
From the results of \cite{arvis} we know that this equation has 3 independent solutions: two even $x(z,\theta)$, 
$y(z,\theta)$ and one odd $\xi(z,\theta)$. Defining $C=y/x$, $\alpha=\xi/x$, we find that $\omega(z,\theta)$ is expressed via 
super Schwarzian derivative, i.e. $w(z,\theta)=\{\alpha; \theta,z\}$ and the consistency conditions on 
$C$ and $\alpha$ are such that $C$ can be represented in terms of $\alpha$ in the following way:
\begin{eqnarray}
C=cA+\gamma A \alpha +\delta\alpha,
\end{eqnarray}
where $A$ is such a function that $(z,\theta)\to (A,\alpha) $ is a superconformal transformation. In a different 
basis $(A, \alpha)$ will be transformed via $SPL_2$ \rf{transff} and hence $(A, \alpha)$ form natural coordinates for a projective structure on $\si$. Therefore we have the following proposition.\\

\noindent{\bf Proposition 2.2} {\it There is a bijection between the set of superprojective structures and the set of superprojective connections on $\si$.}\\

\noindent{\bf 2.3. Connections for vector bundles over super Riemann surfaces.} 
Let us consider a vector bundle $V$ over the super Riemann surface with the fiber $\mathbb{C}^{m|n}_S$. Let $\mathcal{E}^0(\si, V)$ be 
the space of sections on $V$ over $\si$ and let $\mathcal{E}^1(\si, V)$ be the space of 1-form valued sections. As 
usual, the connection is a differential operator
\begin{eqnarray}
d_A(fs)=df\otimes s+(-1)^{|f|}fd_As,
\end{eqnarray}
where $f$ is a smooth even/odd function on $\si$ and $s\in \mathcal{E}^0(\si, V)$. 
Locally, in the chart $(z,\theta)$ the connection has the following form:
\begin{eqnarray}
&&d_A=d+A=d+(\eta A_z+d\theta A_{\theta})+(\bar{\eta}A_{\b z}+d\b\theta A_{\b\theta})=\nonumber\\
&&(\p+\eta A_z+d\theta A_{\theta})+(\b \p+\b \eta A_{\b z}+d\b \theta A_{\b \theta})=\nonumber\\ 
&&(\eta D^A_z+d\theta D_\theta^A)+(\b \eta D^A_{\b z}+d\b \theta D_{\b \theta}^A).
\end{eqnarray}  
We note that we used here the fact that $d=\p+\b \p$ and $\p=\eta\p_z+d\theta D_{\theta}$. The expression for 
the curvature is:
\begin{eqnarray}
&&F=d_A^2=\nonumber\\
&&d\theta d\theta F_{\theta\theta}+\eta d\theta F_{z\theta}+ d\b \theta d\b \theta F_{\b \theta\b \theta}+ \b\eta d\b \theta 
F_{\b z\b \theta}+\nonumber\\
&&\eta\b \eta F_{z\b z}+ \eta d\b \theta F_{z\b \theta}+\b \eta d\theta F_{\b z\theta}+d\theta d\b \theta F_{\theta\b \theta},
\end{eqnarray} 
where $F_{\theta\theta}=- {D^A_{\theta}}^2+D^A_z$, $F_{z\theta}=[D^A_z, D^A_{\theta}]$, $F_{z, \b z}=[D^A_z, D^A_{\b z}]$, $F_{z\b \theta}=
[D^A_z, D^A_{\b \theta}]$, $F_{\theta\b \theta}=-[D^A_{\theta}, D^A_{\b \theta}]$, etc.

It appears that if the connection $d_A$ offers partial flatness, which implies   $F_{\theta\theta}=F_{z\theta}=F_{\b\theta\b \theta}=F_{\b z
\b \theta}=0$, then there is a superholomorphic structure on $V$ (i.e. transition functions of the bundle can be 
made superholomorphic) \cite{RT}. We are interested in the flat superholomorphic connections. In this case, since 
$F_{\theta\theta}=0$, the connection is fully determined by the $D^A_{\theta}$ locally. In other words it is determined by the 
following odd differential operator, which from now on will denote $\nabla$ and call $long$ $superderivative$: 
\begin{eqnarray}\label{fc}
\nabla=D_{\theta}+A_{\theta}(z, \theta), 
\end{eqnarray}
which gives a map: $\mathcal{D}\to End{V}$ so that the transformation properties  for $A_{\theta}$ are: $A_{\theta}\to gA_{\theta}g^{-1}-D_{\theta}g g^{-1}$, where $g$ is a superholomorphic function providing change of trivialization. \\

\noindent{\bf 2.4. $SPL_2$-opers.} In this subsection, we give a description of the first nontrivial superoper. Suppose we have a superprojective structure on $\si$. Naturally we have a structure of a flat $SPL_2$-bundle $\mathcal{F}$ over $\si$, since on on the overlaps there is a constant map to $SPL_2$. Let us study the corresponding flat connection on $\si$. Since $SPL_2$ is a group of superconformal automorphisms of $SC^*$, one can form 
an associated bundle $SC^*_{\mathcal{F}}=\mathcal{F}\times_{SPL_2} SC^*$. This bundle has a global section which is just given by the superprojective coordinate functions $(z, \theta)$ on $\si$. We note that it has nonvanishing (super)derivative at all points. 

One can view $SC^*$ as a flag supermanifold. Namely, consider the group  $SPL_2$ acting in $\mathbb{C}^{2|1}=span(e_1,\xi, e_2)$, where we put the odd vector in the middle. Then $e_1$ is stabilized by the Borel subgroup of upper triangular matrices. Therefore, one can identify $SC^*$ with $SPL_2/B$. Since we have a nozero section of $SC^*_{\mathcal{F}}$, we have a 
$B$-subbundle 
$\mathcal{F}_B$ of a $G$-bundle, where $G$ stands for $SPL_2$. Hence, a superprojective structrure gives the flat $SPL_2$-bundle $\mathcal{F}$ with a reduction $\mathcal{F}_B$. 
However, there is one more piece of data we can use: it is the condition that the (super)derivative of the section of 
$SC^*_{\mathcal{F}}$ is nowhere vanishing. It means that the flat connection on $\mathcal{F}$ does not preserve the 
$B$-reduction anywhere. Let us figure out which conditions does it put on the connection if we choose a trivialization of 
$\mathcal{F}$ induced from the $\mathcal{F}_B$ trivialization.
 As we discussed above, the connection is determined by the following 
odd differential operator:
\begin{eqnarray}
\nabla=D_\theta+
\begin{pmatrix}
 \alpha(z,\theta) & b(z,\theta) &  \beta(z,\theta) \\
 - a(z,\theta) & 0 & b(z,\theta)   \\
  \gamma(z,\theta)   & a(z,\theta)  & -\alpha(z,\theta), 
 \end{pmatrix},
\end{eqnarray}
so that the matrix is in the defining representation of the 
Lie superalgebra of $SPL_2$, namely $osp(2|1)$.  
This operator and its square describe even and odd directions for the tangent vector to $SC^*$. Since we have the condition 
that both of them are nonvanishing, and identifying tangent space with $osp(2|1)/\mathfrak{b}$ (where $\mathfrak{b}$ is the Borel subalgebra), we obtain that $a$ is nonvanishing. It is possible to make $\gamma=0$, 
by redefining $\nabla$ by adding $\mu(\nabla)^2$ with appropriate odd 
function $\mu$, which just corresponds to the choice of superconformal coordinates on $SC^*$. 
We call such a triple $(\mathcal{F}, \mathcal{F}_B, \nabla)$ a $superoper$. We notice that taking the square of the 
odd operator $\nabla$, reducing such even operator 
from $\si$ to the underlying curve $\si^0$ and getting rid of all the odd variables, 
we obtain the oper connection for the $PGL_2$-bundle. Thus superopers can be thought about as ``square roots'' of opers.

Using $B$-valued gauge transformations one can bring $\nabla_{\theta}$ to the canonical form:
\begin{eqnarray}
\nabla=D_\theta+
\begin{pmatrix}
 0 & 0 &  \omega(z,\theta) \\
 -1 & 0 & 0  \\
 0  & 1  & 0 
 \end{pmatrix}.
\end{eqnarray}
Therefore on a superdisc with coordinate $(z,\theta)$ the space of $SPL_2$ superopers can be identified with the space of 
differential operators $D^3_{\theta}-\omega(z,\theta)$. We will see in the next section that the coordinate 
transformations of $\omega$ are the same as in Proposition 2.1.
 
Therefore we see that there is a full analogy with the bosonic case, where the space of $PGL_2$-superopers was identified with the set of projective connections or equivalently with the set of projective structures.

Let us summarize the results of this section in the following theorem.\\

\noindent{\bf Theorem 2.3.} {\it There are bijections between the following three sets on a super Riemann surface $\si$:\\
\noindent i) Superprojective structures\\
\noindent ii) Superprojective connections\\
\noindent iii) $SPL_2$-opers.}\\

\section{Superopers for higher rank superalgebras}

\noindent {\bf 3.1. The definition of superopers.} 
In this section we generalize the results of the previous section to higher rank. Suppose $G$ is a simple algebraic supergroup \cite{berezin}
of adjoint type over Grassmann algebra $S$, $B$ is its Borel subgroup, $N=[B, B]$, so that for the corresponding Lie superalgebras we have ${\mathfrak n}\subset{\mathfrak b}\subset{\mathfrak g}$. 
Note that $\mathfrak{g}=S\otimes \mathfrak{g}^{\rm{red}}$, where 
$\mathfrak{g}^{\rm{red}}$ is a simple Lie superalgebra over $\mathbb{C}$. 
As usual, $H=B/N$ with the Lie algebra $\mathfrak{h}$  
and there is a decomposition: $\mathfrak{g}=\mathfrak{n}_-\oplus\mathfrak{h}\oplus\mathfrak{n}_+$. The corresponding generators of simple roots will be denoted as usual: $e_1, \dots, e_l$; $f_1, \dots, f_l$. We are interested in the superalgebras, which have a pure fermionic system of simple roots, namely $\mathfrak{psl}(n|n)$, $\mathfrak{sl}(n+1|n)$, $\mathfrak{sl}(n|n+1)$, $\mathfrak{osp}(2n\pm 1|2n)$, $\mathfrak{osp}(2n|2n)$, $\mathfrak {osp}(2n+ 2|2n)$ with $n\ge 0$ and $D(2,1;\alpha)$ with $\alpha\neq 0, \pm 1$. Moreover, a necessary ingredient for our construction is the presence of the embeddining of superprincipal 
$osp(1|2)$ subalgebra \cite{FRSp}, \cite{delduc}, namely that for $\chi_{-1}=\sum_i f_i$ and $\check \rho=\sum_i\check{\omega}_i$, where $\check {\omega}_i$ are fundamental coweights, there is such $\chi_{1}$ that makes a triple $(\chi_1, \chi_{-1}, \crho)$ an $osp(1|2)$ superalgebra. Almost all series of 
superalgebras from the list above allow such an embedding, however,  $\mathfrak {psl}(n|n)$ does not and we do not consider these series in this article.  
    
As in the standard bosonic case we define an open orbit 
${\bf O}\subset[\mathfrak{n}, \mathfrak{n}]^{\perp}/\mathfrak{b}$ consisting of vectors, stabilized by $N$ and such that all the negative root components of these vectors with respect 
to the adjoint action of $H$ are non-zero.

Let us consider a principal $G$-bundle $\mathcal{F}$ over $X$, which can be a super Riemann surface $\si$ or a formal superdisc $SD_x={\rm Spec} S[\theta][[z]]$, or a 
punctured superdisc ${D^S_x}^{\times}={\rm Spec} S[\theta] ((z))$ (see e.g. \cite{leites} or \cite{kapranov} for the definitions of the spectra of supercommutative rings), and its reduction $\mathcal{F}_B$ to the Borel subgroup $B$. We assume that it has a flat connection determined by a long 
superderivative $\nabla$ (see \rf{fc}). According to the example, considered in section 2 we do not want $\nabla$ to 
preserve $\mathcal{F}_B$. However, in the higher rank case this is not enough, so we have to specify extra conditions. 
Namely, suppose $\nabla'$ is another long superderivative, which preserves $\mathcal{F}_B$. Then we require that the 
difference $\nabla'-\nabla$ has a structure of superconformal field of dimension $1/2$ with values in the associated bundle 
$\mathfrak{g}_{\mathcal{F}_B}$. We can project it onto  $(\mathfrak{g}/\mathfrak{b})_{\mathcal{F}_B}\otimes \mathcal{D}^{-1}$. Let us   
denote the resulting $(\mathfrak{g}/\mathfrak{b})_{\mathcal{F}_B}$-valued superconformal field as $\nabla/\mathcal{F}_B$. 
Now we are ready to define the superoper, which is a natural generalization of the oper.  

A $G-superoper$ on $X$ is the triple $(\mathcal{F}, \mathcal{F}_B, \mathcal{\nabla})$, where $\mathcal{F}$ is a principle 
$G$-bundle, $\mathcal{F}_B$ is its $B$-reduction and $\nabla$ is a long superderivative on $\mathcal{F}$, such that 
$\nabla/\mathcal{F}_B$ takes values in ${\bf O}_{\mathcal{F}_B}$. 

Locally this means that in the coordinates $(z, \theta)$ and with respect to the 
trivialization of $\mathcal{F}_B$, the structure of the 
long superderivative is:
\begin{eqnarray}\label{sops}
D_{z,\theta}+\sum^l_{i=1}a_i(z, \theta)f_i+\mu(z,\theta),
\end{eqnarray}
where each $a_i(z, \theta)$ is an even nonzero function (meaning that these functions  have nonzero body and are invertible) and 
$\mu(z, \theta)$ is an odd $\mathfrak{b}$-valued function. Therefore locally on the open subset $U$, where we chose 
coordinates $(z, \theta)$, the space of $G$-superopers on $U$, which will be denoted as 
$\sop(U)$, can be characterized the space of all odd 
operators of type $\rf{sops}$ modulo gauge transformations from $B(R)$ group, where $R$ are either analytic or algebraic 
functions on $U$. \\

\noindent {\bf 3.2. Coordinate transformations and other properties.} Let us notice that one can use the $H$-action 
to make the operator \rf{sops} look as follows:
\begin{eqnarray}\label{sopst}
D_{\theta}+\sum^l_{i=1}f_i+\mu(z,\theta),
\end{eqnarray}
where $\mu\in \mathfrak{b}(R)$. Therefore the space $s{\rm Op}_G(U)$ can be considered as the quotient of the space of operators of the form \rf{sopst} 
(denoted as $\widetilde{\rm{sOp}}_G(U)$) 
by the action of $N(R)$. As in the pure bosonic case, $\check{\rho}$ gives a principal gradation 
(for those classes of superalgebras we consider), i.e. we have a direct sum decomposition $\mathfrak{b}=\oplus_{i\ge 0}\mathfrak{b}_i$. Moreover, let us remind that we denoted  
$\chi_{-1}=\sum^l_{i=1} f_i$ and there exists a unique element $\chi_{1}$ of degree $1$ in $\mathfrak{b}$, such that $\chi_{\pm 1}, 
\check{\rho}$ generate $osp(1|2)$ superalgebra. 
Let $\tilde{\chi}_k$ $(k=1,\dots, l)$ (which can be either odd or 
even), so that $\tilde{chi}_2=\chi_1^2$, be the basis of the 
space of the $ad(\chi_1)$ invariants. We note, that the decompositions of $\mathfrak{g}$ with respect to the adjoint action of such 
$osp(1|2)$ triple were studied in \cite{FRSp}. 
Based on that, we have the following Lemma which is proved in a 
similar way as in \cite{DS} (see also Lemma 4.2.2 of \cite{FL}).\\

\noindent{\bf Lemma 3.1.} {\it The gauge action of $N(R)$ on $\widetilde{s\rm{Op}_G}(U)$ is free and 
each gauge equivalence class contains a unique operator of the form \rf{sopst} with 
\begin{eqnarray} 
\mu(\theta, z)=\sum^l_{i=1}g_i(z, \theta)\tilde{\chi}_i,
\end{eqnarray}
where $g_i$ has opposite parity to $\chi_i$.}\\

Now let us discuss the transformation properties of operators $\widetilde\sop(U)$. Assume we have a superconformal coordinate change $(z, \theta)=(f(w,\xi), \alpha(w, \xi))$. 
Then according to the transformations of the long derivative we have
\begin{eqnarray}
&&\nabla=\\
&&D_{\xi}+(D_{\xi}\alpha)(w,\xi)\chi_{-1}+(D_{\xi}\alpha)(w, \xi)(\mu(f((w,\xi),\alpha(w, \xi)).\nonumber
\end{eqnarray}
Considering 1-parameter subgroup ${\mathbb{C}_S^\times}^{1|1}\to H$ which corresponds to $\check{\rho}$, 
applying adjoint transformation with $\check{\rho}(D_{\xi}\alpha)$ we obtain:
\begin{eqnarray}
&&D_{\xi}+\\
&&\chi_{-1}+(D_{\xi}\alpha)(w, \xi)Ad_{\check{\rho}(D_{\xi}\alpha)}\cdot\mu(f((w,\xi), \alpha(w, \xi))-
\frac{\p_w\alpha(w,\xi)}{D_{\xi}\alpha(w, \xi)}\check{\rho}.\nonumber
\end{eqnarray}
This gives us the gluing formula for  superopers on any super Riemann surface $\si$. 

 Consider the $H$-bundle $\mathcal{D}^{-\check\rho}$ on $\si$, which is determined by the  property that the line bundle $\mathcal{D}^{-\crho}\times \mathbb{C}_\lambda$ is $\mathcal{D}^{-\langle\crho, \lambda\rangle}$, where 
$\lambda$ is from the lattice of characters and $\mathbb{C}_{\lambda}$ is the corresponding 1-dimensional representation.

The coordinate transformation formulas for superoper connection immediately lead to another characterization  of this bundle via $\mathcal{F}_B$-reduction. The following statement is the 
supersymmetric version of Lemma 4.2.1 of \cite{FL}. \\ 

\noindent{\bf Lemma 3.2.} {\it The H-bundle $\mathcal{F}_H=\mathcal{F}_B\times_B H=\mathcal{F}_B/N$ is isomorphic to $\mathcal{D}^{-\crho}$.}\\

Now one can derive the transformation properties for the canonical representatives of 
opers from Lemma 3.1, which will provide the transformation formulas for $g_1,\dots, g_n$. In order to do that, one needs to apply to the operator \rf{sopst} the gauge transformation of the form
\begin{eqnarray}\label{tranfun}
\exp\Big({\kappa\chi_1-\frac{1}{2}(D\kappa)[\chi_1, \chi_1]}\Big)\check{\rho}(D_{\xi}\alpha),
\end{eqnarray}
where $\kappa=\frac{\p_w\alpha(w,\xi)}{D_{\xi}\alpha}$.
Then we have that 
\begin{eqnarray}\label{transf}
&&\tilde g_1(w, \xi)= g_1(w, \xi)(D_\xi\alpha)^2,\nonumber\\
&&\tilde g_2(w, \xi)=g_2(w,\xi)(D_\xi\alpha)^3+\{\alpha; w,\xi\},\nonumber\\
&&\tilde g_j(w, \xi)=g_j(w,\xi)(D_\xi\alpha)^{d_j+1}, \quad j>2.
\end{eqnarray}
Therefore \rf{tranfun} are transition functions for $\mathcal{F}_B$ and $\mathcal{F}$  
bundles. \\

\noindent{\bf Remark.} Note that the $g_1$-term is absent in the $osp(1|2)$, however it often appears in the higher rank. The first example is 
$sl(2|1)\cong osp(2|2)$.\\

The formulas \rf{transf} give the following description of the space of superopers:
\begin{eqnarray}
s{\rm{Op}}_G(\si)\cong sProj(\si)\times\oplus^l_{j=1, j\neq 2}\Gamma(\si, \mathcal{D}^{-d_j-1}),
\end{eqnarray} 
where $sProj(\si)$ stands for superprojective connections on $\si$.

In the previous section we indicated that in the $osp(1|2)$ one can introduce the oper related to a superoper, by considering $\nabla^2$, then stripping it from the $\theta$ and $S$ dependence, 
we obtain that the resulting $\overline{\nabla^2}$ has all the needed properties of $sl(2)$ oper on the curve $X$ which is a base manifold for $\si$. 

A similar construction is possible in the higher rank case. Let $^0G$ be 
the reductive group, which is a base manifold for $G$.  Due to the structure of the coordinate transformations we derived above, we find out that indeed $\overline{\nabla^2}=\nabla^2$ defines an oper on $X$. We refer to this object as $^0G$-$oper$, 
$associated$ $with$ $the$ $G$-superoper, which we will 
denote as triple 
$($ ${^0\mathcal{F}},\overline{\nabla^2}$,$ {^0{\mathcal{F}}_B})$, where $^0\mf$,  $^0\mf_B$ denote the appropriate pure even reductions of the principal bundles.

\section{Superopers with regular singularities, Miura superopers and 
Bethe ansatz  equations}

\noindent{\bf 4.1. Superopers with regular singularities}. Consider a point $x$ on on the superc Riemann surface $\si$ and the formal superdisc $SD_x$ around that point with the coordinates $(z,\theta)$. 
We define a $G$-superoper with regular singularity on $SD_x$ as 
an operator of the form 
\begin{eqnarray}\label{sinso}
D_{\theta}+\sum a_i(z, \theta)f_i+\Big(\mu_1(z)+\frac{\theta}{z}\mu_0(z)\Big),
\end{eqnarray} 
modulo the $B(\mathcal{K}_x)$-transformations ($\mathcal{K}_x=\mathbb{C}[\theta]((z))$), where $a_i(z, \theta)\in \mathcal{O}_x$ are nowhere vanishing and invertible, 
$\mu_i(z, \theta)\in \mathfrak{b}(\mathcal{K}_x)$ $(i=0,1)$, such that 
the bodies of $\mu_i(z, \theta)$, i.e. $\overline{\mu_i}
\in \mathfrak{b}^{\rm{red}}(\mathcal{O}^{\rm{red}}_x)$. 
As before, one can eliminate $a_i$-dependence via $H$-transformations, therefore we can talk about 
$N(\mathcal{K}_x)$ equivalence class of operators of the type \rf{sinso} with $a_i=1$. Let us denote by $sOp_G^{RS}(SD_x)$ the space of superopers with regular singularity. Clearly, we have the embedding: $sOp^{RS}_G(SD_x)\subset sOp_G(SD^\times_x)$.


The $^0G$-oper, corresponding to $G$-superoper \rf{sinso} is the oper with regular singularity. It has the following form:
\begin{eqnarray}
\p+\chi^2_{-1}+[\chi_{-1},\overline{\mu_1}]+ (\overline{\mu_1})^2+\frac{1}{z}(\overline{\mu_0}),
\end{eqnarray}
which can be transformed to the standard form via 
the gauge transformation by means of $\frac{\rho}{2}(z)$:
\begin{eqnarray}
\p_z+\frac{1}{z}\Big(\chi^2_{-1}-\frac{\check\rho}{2}+Ad_{\frac{\check{\rho}}{2}(z)}\cdot\bar{\mu}_0\Big)+v(z),
\end{eqnarray}
where $v(z)$ is regular.


Denoting $-\check\lambda$ the projection    
of $\mu_0$ on $\mathfrak{h}$, we find that the residue of this 
differential operator is equal to $\chi_{-1}^2-\lambda-\2\check\rho$, however since this is an oper, only the corresponding class in $\mathfrak{h}/W$ is well defined, and we denote it as $(-\lambda-\2\check\rho)_W$, i.e. this oper belongs to ${\rm Op}^{RS}_G(D_x)_{\check\lambda}$, see e.g. \cite{Fe}. 

Let us refer to the space of superopers with regular singularity such that  $\bar{\mu_0}(0)=\check\lambda$, as $s{\rm Op}_G^{RS}(D_x)_{\lambda}$. 

If we consider the representation $V$ of $G$ one can talk about a system of differential equations $\nabla\cdot \phi_V(z, \theta)$ and their monodromy like in the pure even case. 

Let $\check{\lambda}$ be the dominant integral coweight and let us introduce the following class of operators:
\begin{eqnarray}\label{lsop}
\nabla=D_{\theta}+\Big(\sum a_i(z, \theta)f_i+\mu(z, \theta)\Big),
\end{eqnarray} 
where $a_i=z^{\langle \alpha_i, \check{\lambda}\rangle}(r_i(\theta)+z(\dots))$, so that the body of $r_i$ is nonzero, and $\mu(z, \theta)\in \mathfrak{b}(\mathcal{O}_x)$.  
We call the quotient of the space of operators above by 
the action of $B(\mathcal{O}_x)$ as 
$\sop (SD_x)_{\lambda}$. 

The following Lemma is an analogue of Lemma 2.4. of \cite{Fe}.\\

\noindent{\bf Lemma 4.1.} {\it There is an injective map i
: $\sop(SD_x)_{\check\lambda}\to sOp(SD^{\times}_x)$, so that $\rm{Im}$i $ \subset \sop^{RS}(SD_x)_{\check\lambda}$.  The image of i is a subset in the set of those elements of $\sop^{RS}(SD_x)_{\check\lambda}$, such that the resulting oper has a trivial monodromy around x.}\\

\noindent{\bf Remark.} Notice that the superopers  corresponding to $s{\rm Op}_G(SD_x)_{\check\lambda}$ belong to 
${\rm Op}_G(D_x)_{\check\lambda}$. However, here $\check\lambda$ is the integral dominant weight for Lie superalgebra. If we consider $\lambda$ to be an integral dominant weight for the underlying Lie algebra, the monodromy for the corresponding superoper would not be necessarily trivial: the expression will include the half-integer powers of $z$ and the monodromy will correspond to the reflection: $\theta\to -\theta$.\\

\noindent{\bf 4.2. Miura superopers.} Miura superoper is defined in complete analogy with the pure even case. Namely, {\it Miura G-superoper} 
is a quadruple $(\mathcal{F}, \nabla, \mathcal{F}_B, \mathcal{F}'_B)$ 
where the triple $(\mathcal{F}, \nabla, \mathcal{F}_B)$ is a $G$-superoper 
and ${\mathcal{F}}'_B$ is another B-reduction preserved by $\nabla$. Let us denote the space of such superopers as $\smop(\si)$. 

Such $B$-reductions of $\mathcal{F}$ are completely determined by the B-reduction of the fiber $\mathcal{F}_x$ at any point $x$ on $\si$ and a set of all such reductions is given by $(G/B)_{\mathcal{F}_x}=\mathcal{F}_x\times_G G/B=(G/B)_{\mathcal{F}'_x}
$. Then if superoper $\xi$ has the regular singularity and a trivial 
monodromy, then there is an isomorphism between the space of Miura opers for such $\xi$ and $(G/B)_{\mathcal{F}'_x}$.

The structure of the flag manifold $G/B$ is usually quite complicated \cite{penkov},\cite{manin}, however we just need the structure determined 
by its "body", i.e. $^0G/^0B$. For the pure even flag variety 
$^0G/$ $^0B$, we have the standard Schubert cell decomposition, 
where cells $^0S_w$=$^0Bw_0$$w$ $^0B$ are labeled by the Weyl group elements 
$w\in W$  and $w_0$ is the longest element of the Weyl group (from now on when we say Weyl group, we mean only the Weyl group corresponding to pure even Weyl reflections of the $^0G$ root system).  

Let us denote $S_w$ the preimage of 
$P:G/B$$\to$ $^0G/$ ${^0B}$. 
We assume that the preimage of a big cell $^0Bw_0^0B$ allows factorization $Bw_0B$.    
The B-reduction $\mathcal{F}'_B$ defines a point in $G/B$.
We say that B-reductions $\mathcal{F}_{B,x}$ and $\mathcal{F}'_{B,x}$ are in relative position $w$ if $\mathcal{F}_{B,x}$ belongs to $\mathcal{F'}_x\times_B S_w$. When $w$=1, we say that $\mf_x$, $\mf_x'$ are in generic position. A Miura superoper is called 
generic at a given point $x\in \si$ if the $B$-reductions $\mf_{B,x}$, $\mf'_{B,x}$. are generic. Notice that if a Miura superoper is generic at $x$, it is generic in the neighborhood of x. 
We denote the space of Miura superopers on U as $\smop(U)_{gen}$. It is clear that the reduction of Miura superoper to $(^0\mathcal{F}, \overline{\nabla^2}, ^0\mathcal{F}_B, ^0\mathcal{F}'_B)$ gives a Miura oper. 

Therefore the following Proposition holds, which follows directly from the reduction to the pure even case, although one can also go along the lines of the proof of Lemma 2.6. and Lemma 2.7 of \cite{Fe}.\\

\noindent{\bf Proposition 4.2.} {\it i) The restriction of the Miura superoper to the punctured disk is generic.\\
ii) For a generic Miura superoper $(\mathcal{F}, \nabla, \mathcal{F}_B, \mathcal{F}'_B)$ the $H$-bundle $\mathcal{F}'_H$ is isomorphic to 
$w_0^*(\mathcal{F}_h)$
}\\ 

As in the even case we can define an $H$-connection associated to Miura superoper on $\mathcal{F}_H\cong \mathcal{D}^{-\check\rho}$, which 
is determined by $\tilde\nabla=D_{\theta}+u(z,\theta)$, where $u$ is $\mathfrak{h}$-valued function. Under the change of coordinates $(z, \theta)=(f(w,\xi), 
\alpha(w,\xi))$, the long superderivative transforms as follows:
\begin{eqnarray}
D_{\xi}+D_{\xi}\alpha \cdot u(f(w, \xi), \theta(w,\xi))-
\frac{\p_w\alpha(w,\xi)}{D_{\xi}\alpha)}\check{\rho}.
\end{eqnarray} 
Let us call the resulting morphism $\smop(U)_{gen}$ to the space $Conn_U$ of the described above flat $H$-connections on $U$ as $\mathbf{a}$. 
Now suppose we are given a long superderivative $\tilde\nabla$ on $H$-bundle $\mathcal{D}^{-\check\rho}$, one can construct a generic superoper as follows. Let us set  
$\mathcal{F}=\mathcal{D}^{-\check \rho}\times_H G$, $\mf_B=\mathcal{D}^{-\check \rho}\times_H B$. Then, defining $\mf_B'$ as 
$\mathcal{D}^{-\check \rho}\times_H w_0 B$ and the long
 superderivative on $\mathcal{F}$ as $\nabla=\chi_{-1}+\tilde\nabla$,  
 we see that the constructed quadruple $(\mf, \nabla, \mf_B, \mf_B')$ is a generic Miura oper.

Therefore, we obtained the following statement which is analogue  
of Proposition 2.8 of \cite{Fe}.\\

\noindent{\bf Proposition 4.3.} {\it The morphism ${\bf a}: \smop(U)_{gen}\to Conn_U$ is an isomorphism of algebraic supervarieties.}\\

Similarly one can define the space of Miura $G$-superopers of coweight $\check\lambda$ on $SD_x$ via the same definition applied to $\sop(SD_x)_{\check\lambda}$. Again, we have isomorphism $\smop(SD_x)_{\lambda}\cong\sop(SD_x)_{\lambda}\times (G/B)_{\mathcal{F}_x'}$. 
We define relative positions as in the case of standard Miura superopers (
$\check\lambda=0$) and let $\smop(SD_x)_{\check\lambda, gen}$ denote the variety of generic Miura opers of weight $\lambda$. 

Finally, there is an analogue of Proposition 4.3 in this case. 
Let $Conn^{RS}_{SD_x, \check\lambda}$ denote the set of  of long derivatives on the H-bundle $\mathcal{D}^{-\check\rho}$ with regular singularity and residue $-\check\lambda$, namely the long derivatives of the form:
\begin{eqnarray}
\tilde \nabla=D_{\theta}+\frac{\theta}{z}\check\lambda+u(z,\theta),
\end{eqnarray} 
where $u(z,\theta)\in \mathfrak{h}[[z,\theta]]$. Then as before, one can construct a connection $\nabla=\tilde\nabla+\chi_{-1}$ and making the gauge transformation with $\check\lambda(z)$ we obtain the connection from $\sop(SD_x)_{\check\lambda}$. Therefore, there is an isomorphism between $Conn^{RS}_{SD_x, \lambda}$ and $\smop(SD_x)_{gen,\check\lambda}$.\\

\noindent{\bf 4.3. Miura superopers with regular singularities on $SC^*$.}
 First, let us consider a Miura superoper of coweight $\hat\lambda$ on the disc $SD_x$. 
 Assume, it is not generic, but $\mf'_{B,x}$ has the relative position $w$ with $\mf_{B,x}$ at $x$. Let us denote the space of all such Miura superopers by $\smop{(SD_x)}_{\check\lambda,w}$.
 
 From previous subsection we know that each such Miura superoper corresponds to some $H$-connection on $\mathcal{D}^{-\check\rho}$ over $SD_x^{\times}$. Using the results from the pure even case, one can show that the corresponding $H$-connection has the form 
\begin{eqnarray}\label{nu}
D_{\theta}+\frac{\theta}{z}\check\nu+f(z,\theta) ,
\end{eqnarray}
 where $\check\nu-\2\check\rho=w(-\check\lambda-\2\check\rho)$,  $w$ defines the relative position at $x$,  $f(u,\theta)$ is such that the body of its superderivative is regular in $z$, i.e. $\overline{D_{\theta}f(z,\theta)}\in \mathfrak{h}^{\rm red}[[z]]$. 
 Let us call the space of such connections by $Conn^{RS}_{SD_x, \check\lambda, w}$.

Therefore, we can construct a map ${\bf b}^{RS}_{\lambda,w}: Conn^{RS}_{SD_x, \check\lambda, w}\to \sop^{RS}(SD_x)$
 similarly to the previous subsection, by constructing  the triple 
 $(\mf, \nabla, \mf_B)$ via identification 
 $\mf=\mathcal{D}^{-\crho}\times_H G$, 
 $\mf_B=\mathcal{D}^{-\crho}\times_H B$ and 
 $\nabla=\t\nabla+\chi_{-1}$, where $\t \nabla\in Conn^{RS}_{SD_x, \check\lambda, w}$. We denote by 
$Conn^{reg}_{SD_x,\check\lambda,w}$ the preimage of $\sop{(SD_x)}_{\check\lambda,w}$ under this morphism, therefore we have the map:
${\bf b}_{\lambda,w}: 
Conn^{reg}_{SD_x, \check{\lambda}, w}\to \smop^{RS}(SD_x)_{\check{\lambda}}$, so that in the quadruple $(\mf, \nabla, \mf_B, \mf'_B)$ first three terms are as above and $\mf'_B=\mathcal{D}^{-\crho}\times_H w_0B$. If we denote 
$\smop^{RS}(SD_x)_{\check{\lambda}, w}$ those Miura superopers of coweight $\check{\lambda}$ which have the relative position $w$ at $x$, then the following Proposition is true, based on the results from the pure even case (see Proposition 2.9 of \cite{Fe}).\\

\noindent{\bf Proposition 4.4.} {\it For each $w\in W$, ${\bf b}_{\cla, w}$ is an isomorphism of supervarieties $Conn^{reg}_{SD_x,\check\lambda,w}$ and $\smop{(SD_x)}_{\check\lambda,w}$}\\

Let us now consider the case of $\check\lambda=0$ and assume that the relative position is  given by $s_{2\alpha_i}$, where $\alpha_i$ is a simple black root.  In local coordinates, the corresponding $H$-connection will be given by the differential operator:
\begin{eqnarray}\label{simple}
\tilde \nabla=D_{\theta}+\frac{\theta}{2z}\check{\alpha}_i+u(z,\theta),
\end{eqnarray}
where $u(z,\theta) \in \mathfrak{h}[\theta]((z))$ and $u(z,\theta)=u_1(z)+\theta u_0(z)$ and $\overline{u_0}(z)\in \mathfrak{h}^{\rm red}[[z]]$.  
Then applying the gauge transformation   
\begin{eqnarray}
\exp\Big(-\frac{\theta}{2z}e_{i}+\frac{1}{4z}e^2_{i}\Big)
\end{eqnarray}
to the Miura superoper $\tilde \nabla+\chi_{-1}$, we obtain that the resulting element 
of $\sop{(SD_x)}_{\check\lambda,s_{2\alpha_i}}$ gives the element ${\rm Op}_G{D_x}_{\check\lambda,s_{2\alpha_i}}$ if $\langle\check\alpha_{i}, \overline{u}_0(0)\rangle=0$. 
If we consider the associate bundle corresponding to the 3-dimensional representation of the $osp(2|1)$ triple $\{e_i, f_i, \check\alpha_i\}$ , writing explicitly all the solutions we find that this condition is also a necessary one. Namely, the following Proposition holds.\\

\noindent{\bf Proposition 4.5.} {\it A superoper corresponding to the H-connection given by 
\rf{simple} corresponds to 
${\rm Op}_G{(D_x)}_{\check\lambda,s_{2\alpha_i}}$ if and only if $\langle\check\alpha_{i}, \overline{u_0}(0)\rangle=0$.}\\

Now we are ready to study superopers with regular singularities over the super Riemann surface $SC^*$.  Let us consider $\mathcal{Z}_1=(z_1,\theta_1),\dots,\mathcal{Z}_N=(z_N, \theta_N)$ on $SC^*$. Also, 
let   $\check\lambda_1,\dots, \check\lambda_N,\check\lambda_{\infty}$ be the set of dominant coweights of $\mathfrak{g}$.
Let us consider the H-connections on $SC^*$ with regular singularities at the points $\mathcal{Z}_1, \dots, \mathcal{Z}_N, (\infty,0)$  and a finite number of other points 
$\sw_1=(w_1, \xi_1), \dots, \sw_n=(w_m,\xi_m)$ such that the residues of the corresponding even $H$-connection at $z_i$, $w_j$, $\infty$ are equal to 
$-y_i(\cla+\2\crho)+\crho$, $-y'_j(\crho)+\2\crho$, 
$-y_i(\cla_{\infty}+\2\crho)+\2\crho$, where $y_i, y_j', y_{\infty}\in W$. In other words, we are considering the 
H-connections determined by the differential operator of the following type:
\begin{eqnarray}
&&D_{\theta}-\Big(\sum^N_{i=1}\frac{\theta-\theta_i}{z-z_i+\theta\theta_i} (y_i({\cla}+\frac{\crho}{2})-\frac{\crho}{2})+\nonumber\\
&&\sum^m_{j=1}\frac{\theta-\xi_j}{z-w_j+\theta\theta_j} (y'_j(\frac{\crho}{2})-\frac{\crho}{2}\Big) +\rm{nilp} 
\end{eqnarray}
on $SC^*{\backslash}\infty$, where 
{\rm nilp} stands for elements $f(z,\theta)$ from 
$\mathfrak{h}[\theta]((z))$ such that $\overline{f(z,\theta)}=
\overline{D_{\theta}f(z,\theta)}=0$. 
Let us study its behaviour at infinity. Any connection $D_{\theta}+\alpha(\theta,u)$ on $\mathcal{D}^{-\crho}$ has the following expansion with respect to the coordinates $(u,\eta)=(\frac{-1}{z}, \frac{\theta}{z})$:
\begin{eqnarray}
 D_{\eta}+u^{-1}\alpha(-u^{-1}, -\eta u^{-1})+u^{-1}{\eta}\crho.
\end{eqnarray}
Therefore, considering $\frac{\eta}{u}$-coefficient in the expansion, we obtain the following constraint:
\begin{eqnarray}\label{constraint}
&&\sum^{N}_{i=1}(y_i(\cla+\frac{\crho}{2})-\frac{\crho}{2})+
\sum^{m}_{i=1}(y'_i(\frac{\crho}{2})-\frac{\crho}{2})=\nonumber\\
&&y'_{\infty}(-w_0(\cla_{\infty})+\frac{\check{\rho}}{2})-\frac{\check{\rho}}{2},
\end{eqnarray}
where $y'_{\infty}w_0=y_{\infty}$.  
This expression is expected from the consideration of the pure even case \cite{Fe}. 

Let us denote the set of the considered above $H$-connections by  $Conn(SC^*)^{RS}_{(\sz_i), (\infty,0); \check\lambda_i, \check{\lambda}_{\infty}}$.

Now one can associate to any such connection a $G$-oper on $SC^*$ with regular singularities at the points $(\sz_i)$, 
$(\sw_j)$, $(\infty,0)$ by setting, in familar way,  $\mathcal{F}=\mathcal{D}^{-\crho}\times_H G$,  $\mathcal{F}_B=\mathcal{D}^{-\crho}\times_H B$. 

Let us denote the set of superopers with regular singularities at $\sz_1\dots \sz_N, (\infty,0)$, whose restriction to the formal superdisc at any point $\sz_i$ or $(\infty,0)$ belongs to the space $\sop(SD_{\sz_i})_{\check \lambda}$ or 
$\sop(SD_{(\infty,0)})_{\cla_{\infty}}$, by 
$\sop(SC^*)_{(\sz_i), (\infty,0); \check\lambda_i, \check{\lambda}_{\infty}}$. 

Then let $Conn(SC^*)_{(\sz_i), (\infty,0); \check\lambda_i, \check{\lambda}_{\infty}}\subset Conn(SC^*)^{RS}_{(\sz_i), (\infty,0); \check\lambda_i, \check{\lambda}_{\infty}}$ be those $H$-connections with regular singularities, which are associated to \\
$\sop(SC^*)_{(\sz_i), (\infty,0); \check\lambda_i, \check{\lambda}_{\infty}}$ under the above correspondence. 
Therefore we have the map
\begin{equation}
Conn(SC^*)_{(\sz_i), (\infty,0); \check\lambda_i, \check{\lambda}_{\infty}}\to
\sop(SC^*)_{(\sz_i), (\infty,0); \check\lambda_i, \check{\lambda}_{\infty}}.
\end{equation}
We can construct a Miura superoper associated with the image of this map, namely $\mf'_B=\mathcal{D}^{-\check{\rho}}\times_H w_0B$. 
Therefore, this map can be lifted to 
\begin{eqnarray}
&&{\bf b}_{(\sz_i), (\infty,0); \check\lambda_i, \check{\lambda}_{\infty}}:\nonumber\\
&&Conn(SC^*)_{(\sz_i), (\infty,0); \check\lambda_i, \check{\lambda}_{\infty}}\to
\smop(SC^*)_{(\sz_i), (\infty,0); \check\lambda_i, \check{\lambda}_{\infty}}.
\end{eqnarray}
Similarly to the pure even case, one can argue that this map is an isomorphism. 
Notice that for a given superoper $\tau\in \sop(SC^*)_{(\sz_i), (\infty,0); \check\lambda_i, \check{\lambda}_{\infty}}$ (because of the absence of nontrivial monodromy), the 
space $\smop(SC^*)_{\tau}$ 
of the corresponding Miura superopers is isomorphic to $G/B$.

Similarly to the argument in the pure even case, we obtain the following theorem, which is an analogue of Theorem 3.1 of \cite{Fe}.\\

\noindent{\bf Theorem  4.6.}{\it The set of all connections  $Conn(SC^*)_{(\sz_i), (\infty,0); \check\lambda_i, \check{\lambda}_{\infty}}$, which correspond to a given oper $\tau\in \sop(SC^*)_{(\sz_i), (\infty,0); \check\lambda_i, \check{\lambda}_{\infty}}$, is isomorphic to the set of points of the flag variety G/B.}\\

\noindent{\bf 4.4. $SPL_2$-superopers and super Bethe ansatz equations.}
In this section we return back to the simplest nontrivial example of the superoper, related to supergroup $SPL_2$.  In the previous section we obtained that for a fixed superoper $\tau$ one can trivialize $\mathcal{F}$ by using the fiber at $(\infty, 0)$. Therefore we have the trivialization of $G/B$- bundle and the map: $\phi_{\tau}:SC^*\to G/B$, so that $(\infty,0)$ maps into the point orbit of $G/B$. Also, in the case $G=SPL_2$,  $G/B\cong SC^*$. 

Similar to the pure even case, let us call the superoper $\tau$ $non-degenerate$ if 
i) $\phi_{\tau}(\mathcal{Z}_i)$ 
is in generic position with 
$B$, 
for any 
$i=1,\dots, N$,  
ii) The relative position of $\phi_{\tau}(x)$ and $B$ is either generic or corresponds to a reflection for all $x\in SC^*\backslash {(\infty,0)}$. Since $PGL_2$ opers are non-degenerate for the generic choice of $z_i$, and those are the opers corresponding to $SPL_2$-superopers, then any $\tau\in s{\rm Op}_{SPL(2)}(SC^*)_{(\sz_i), (\infty,0); \check\lambda_i, \check{\lambda}_{\infty}}$ for the generic choice of $\sz_i$ is non-degenerate. 
Also, let us consider the unique Miura superoper structure for $\tau$, such 
that $\mf_{B, (\infty,0)}$ and $\mf'_{B, (\infty,0)}$ coincide, i.e. correspond to the point orbit in $G/B$.

The corresponding $H$-connections will have the following form:
\begin{eqnarray}
\tilde{\nabla}=D_{\theta}-\sum^N_{i=1}\frac{\theta-\theta_i}{z-z_i-\theta\theta_i} {\cla}_i+
\sum^m_{j=1}\frac{\theta-\xi_j}{z-w_j-\theta\xi_j} \frac{\check{\alpha}}{2}+\rm{nilp},  
\end{eqnarray}
where $\check{\lambda}_i=l_i \check\omega$, so that $l_i\in\mathbb{Z}_+$. 
Imposing the constraint from Proposition 4.5, we obtain that the following equations should hold for the corresponding oper to be monodromy free:
\begin{eqnarray}\label{gbae}
\sum^N_{i=1}\frac{2l_i}
{w_j-z_i}-
\sum^m_{s=1}
\frac{2}{w_j-w_s}=0
\end{eqnarray}
Also, let us recall that the coweights $\cla_i$ should also satisfy \rf{constraint}, which in our case simpifies to:
\begin{eqnarray}
\sum^N_{i=1}l_i-m=l_{\infty}
\end{eqnarray}
Note, that the corresponding $PGL_2$-oper coweights, i.e. $2l_i$ are even: superopers associated with the odd weights will have a monodromy which will correspond to a reflection in $\theta$ variable, as it was explained above. 
The equations \rf{gbae} are exactly the Bethe ansatz equations for $osp(2|1)$ Gaudin model studied in \cite{kulish}.\\

\section{Some remarks}
In this article, we studied superopers for superalgebras with pure odd simple 
root system. However, one can define a similar object for other types of superalgebras, just in such case it can be only locally defined (i.e. on a superdisc). The analogue of the expression \rf{sops} will be:
\begin{eqnarray}
\nabla=D_{z,\theta}+\sum_{e}a_e(z, \theta)f_e+\sum_{o} \theta a_o(z, \theta)f_o+
\mu(z,\theta),
\end{eqnarray}
where the summation is over even and odd roots correspondingly and 
$a_{e,f}(z,\theta)$ are the even functions of $z, \theta$ with nonzero body. 
The resulting connection cannot be defined globally on the super Riemann surface, however the operator $\nabla^2\vline_{\theta=0}$ can give rise to a connection for a G-bundle over a smooth curve underlying the super Riemann surface, while 
$\overline{\nabla^2}$ will give an oper for the underlying semisimple supergroup. 
This construction gives a generalization of opers in the case of any simple superalgebra. 

In this paper we briefly considered an important relation between the spectrum of the Gaudin model and superopers on $SC^*$, which in fact could 
give an example of geometric Langlands correspondence in the case of  superalgebras. For $SPL_2$-superopers and Gaudin model for $osp(2|1)$ the spectrum was determined in fact by the underlying $PGL_2$-oper. Unfortunately so far Gaudin models were not studied in the case of other superalgebras yet, so it is not clear whether such a relation holds for higher rank superalgebras.

We will address these and other important questions in the forthcoming publications.


\begin{thebibliography}{10}
\bibitem{arvis}J.F. Arvis, {\it Classical dynamics of supersymmetric Liouville theory}, Nucl. Phys. {\bf B} 212 (1983) 151-172.
\bibitem{BD} A. Beilinson, V. Drinfeld, {\it Opers}, arXiv: math.AG/0501398.
\bibitem{berezin}F.A. Berezin, {\it Introduction to Superanalysis}, Springer (1987). 
\bibitem{br}M. J. Bergvelt and J. M. Rabin. {\it Supercurves, their Jacobians, and super KP equations}, 
Duke Math. Journal, 98(1), 1999.
\bibitem{crane} L. Crane, J.M. Rabin, {Super Riemann Surfaces: Uniformization and Teichmueller Theory}, Comm. Math. Phys. 113 (1988) 601-623.
\bibitem{delduc} F. Delduc, E. Ragoucy, P. Sorba, {\it Super-Toda Theories and W-algebras from Superspace Wess-Zumino-Witten Models}, Commun. Math. Phys. {\bf 146} (1992) 403-426.
\bibitem{DG} F. Delduc, A. Gallot, {\it Supersymmetric Drinfeld-Sokolov reduction},arXiv: solv-int/9802013.
\bibitem{DS} V. Drinfeld, V. Sokolov, {\it Lie Algebras and KdV type equations}, J. Sov. Math. {\bf 30} (1985) 1975-2036.
\bibitem{FFR} B. Feigin, E. Frenkel, N. Reshetikhin, {\it Gaudin model, Bethe Ansatz and critical level}, Comm. Math. Phys. {\bf 166} (1994) 27-62.
\bibitem{FRSp} L. Frappat, E. Ragoucy, P. Sorba, 
{\it W-algebras and superalgebras from constrained WZW models: a group theeoretical classification}, arXiv: hep-th/9207102. 
\bibitem{Fe} E. Frenkel, {\it Affine Algebras, Langlands Duality and Bethe ansatz}, in Proceedings of International Congress of Mathematical Physics, 
Paris, 1994, International Press, 606-642.
\bibitem{FL} E. Frenkel, {\it Langlands correspondence for loop groups}, CUP, 2007.
\bibitem{Fb} E. Frenkel, {\it Opers on the Projective Line, Flag Manifolds and Bethe Ansatz}, arXiv: math/0308269.
\bibitem{inami} T. Inami, H. Kanno, {\it Lie Superalgebraic Approach to Super Toda Lattice and Generalized Super KdV Equations}, Commun. Math. Phys. {\bf 136} (1991) 519-542.
\bibitem{kapranov} M. Kapranov, E. Vasserot, {\it Supersymmetry and the formal loop space}, arXiv:1005.4466.
\bibitem{kulish} P.P. Kulish, N. Manojlovic,  {\it Bethe vectors of the $osp(1|2)$ Gaudin model },  Lett.Math.Phys. 55 (2001) 77-95. 
\bibitem{kz} P.P. Kulish, A.M. Zeitlin, {\it Group Theoretical Structure and Inverse Scattering Method for super-KdV Equation}, J. Math. Sci {\bf 125}  (2005)203-214.
\bibitem{leites} D. A. Leites, {\it Theory of supermanifolds}, KF Akad. Mauk SSSR, Petrozavodsk (1983).
\bibitem{manin}  Yu.I. Manin, A.A. Voronov,  {\it Supercellular partitions of flag superspaces}, Itogi Nauki i Tekhniki. Ser. Sovrem. Probl. Mat. Nov. Dostizh., 32, VINITI, Moscow, 1988, 27�70. 
\bibitem{penkov} I.B. Penkov, {\it Borel-Weil-Bott theory for classical Lie supergroups}, Itogi Nauki i Tekhniki. Ser. Sovrem. Probl. Mat. Nov. Dostizh., 32, VINITI, Moscow, 1988, 71�124. 
\bibitem{mathieu}P. Mathieu, {\it Super Miura transformations, Super Schwarzian derivatives and Super Hill Operators}, 
in: Integrable and Superintegrable systems, World Scientific (1991) 352-388 
\bibitem{RT} M. Rakowski, G. Thompson, 
{\it Connections on Vector Bundles over Super Riemann Surfaces}, Phys. Lett. {\bf B} 220 (1989) 557-561.
\bibitem{wang} S.-J. Cheng, W. Wang, {\it Dualities and Representations of Lie Superalgebras},  Graduate Studies in Mathematics 144, AMS, 2012
\bibitem{W} E. Witten, {\it Khovanov Homology and Gauge Theory},  arXiv: 1108.3103
\bibitem{Wl} E. Witten, {\it Notes on Super Riemann Surfaces and their Moduli}, arXiv:1209.2459    
\end{thebibliography}
\end{document}